# Optimal Operation of PV-Battery-Diesel MicroGrid for Industrial Loads Under Grid Blackouts

Mansour Alramlawi, *Student Member, IEEE*, Aouss Gabash, *Member, IEEE*,
Erfan Mohagheghi, *Member, IEEE*, Pu Li
Department of Simulation and Optimal Processes
Institute of Automation and Systems Engineering
Ilmenau University of technology
Ilmenau, Germany
mansour.alramlawi@tu-ilmenau.de, aouss.gabash@tu-ilmenau.de, erfan.mohagheghi@tu-ilmenau.de, pu.li@tu-ilmenau.de

*Abstract*—In this paper, we propose an approach to determine the optimal operation strategies for a PV-diesel-battery microgrid covering industrial loads under grid blackouts. A special property of the industrial loads is that they have low power factors. Therefore, the reactive power consumption of the load cannot be neglected. In this study, a novel model of a PV-battery-diesel microgrid is developed considering the active as well reactive power of the microgrid components. Furthermore, an optimization approach is proposed to optimize the active as well reactive power flow in the microgrid for covering the load demand while decreasing the power consumption from the grid, minimizing the diesel generator (DG) operation cost as well as maximizing the consumed power from the PV-array. It has been found that the proposed operation strategy induces a huge reduction of the consumed energy cost and the PV curtailment.

## I. Introduction

Grid blackout is a major problem for industrial sectors in many countries through the world [1]. DGs were found as a proper solution to cover the load demand during blackouts periods [2], [3]. However, DGs have many disadvantages such as a high operational and maintenance costs and environmental pollution. Due to the technology development in the recent years, it has become possible to integrate renewable energy sources with conventional sources to form a local microgrid that is able to provide an uninterruptible power supply [4].

Furthermore, the authors in [5] proved that an optimal sizing of a grid-connected PV-battery system makes it able to provide a reliable power source for households, suffering from intermittent grid. In addition, the work in [6] presented a design approach for sizing PV-battery-diesel microgrid components to provide a reliable energy source considering planned grid outages.

From another perspective, a predictive energy management and control system was proposed in [7] for a PV-battery microgrid connected to an unreliable public grid which suffering from periodic blackouts. The proposed control strategy aimed to minimize the PV power curtailment and increase grid reliability. In [8], a two-stage model predictive controller (MPC) is introduced to optimize the operation of an islanded PV-diesel-battery microgrid and guarantee an uninterrupted power supply for the load. The authors in [9] proposed continuous and ON/OFF control methods to reduce the operation cost of a stand-alone PV-battery-diesel system considering only the DG fuel consumption cost. The obtained results showed that the continuous control strategy is more efficient than the ON/OFF control strategy in reduction of the operation cost. A battery management system for a PV-battery diesel microgrid is proposed in [10] to reduce DG operating hours, controlling battery bank charge and discharge processes and minimizing PV power output fluctuations. In [11], six operation modes were proposed to operate a PV-battery system taking into acount the grid scheduled blackouts problem based on the available power sources priorities. These modes were seasonally selected to decrease the total cost of the consumed energy from the grid using pareto optomization. Recently, an economic MPC was used in [12] to create an optimal power dispatch framework for a PV-battery system taking into account battery lifetime and grid blackouts problem. It is worth to mention that the reactive power has not been considered in the previous studies [7]–[12]. A special property of industrial loads is that they have low power factors. Therefore, the active and reactive power flow within the microgrid should be considered. This leads to difficulties in satisfying the load active and reactive power demand while holding all technical constraints of the microgrid components. Moreover, [13] and [14] have demonstrated that using reactive power capability of renewable energy sources could lead to reduction of the imported reactive power from the grid, decrease the power curtailment from the connected wind turbines and increase system profitability.

Comparing to the above studies, the contributions of this paper can be summarized as follows:

1) A novel model of PV-battery-diesel microgrid is developed to consider both active and reactive power flow and grid blackout problem.

This research is financed by the German Academic Exchange Service (DAAD).

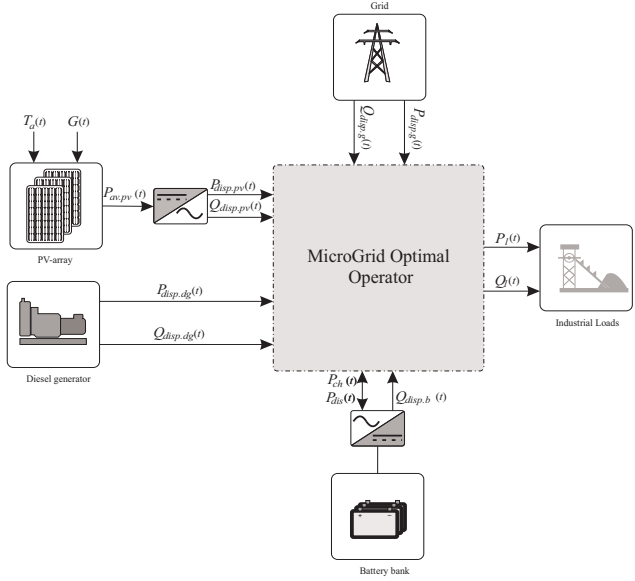

Fig. 1. The proposed PV-battery-diesel microgrid block diagram

2) An optimal operation strategy is introduced to optimize the PV-battery-diesel microgrid that simultaneously minimizes the power consumption cost and maximizes the consumed power from the PV-array.
3) The performance of the proposed model and operation strategy is analyzed at each season to confirm the effectiveness of the developed approach as a potential solution for the grid blackout problem.

The remainder of this paper is formulated as follows. The PV-battery-diesel microgrid model is introduced in section II. In section III optimal operation strategy is described. The results presented and discussed In Section IV. Finally section V concludes with a summary and future work .

## II. PV-BATTERY-DIESEL MICROGRID MODEL

The proposed PV-battery-diesel microgrid block diagram is shown in Fig. 1. It is a grid-connected microgrid consisting of a PV-array with a PV-inverter, a battery bank with a battery inverter, a diesel generator and industrial load.

### A. PV-Array

The PV-cell in the core components of the PV-array, which is reposible to converte solar energy to electrical energy. In this study, the single diode PV-cell model [15] is utilized to determine the maximum output power generated from each PV-cell $P_{max.c}(t)$ in the PV-array. The output power of the PV-cell at each time step is highly depends on the ambient temperature $T_a(t)$ and the global solar irradiance $G(t)$ values. It can be calculated as follows [12]

$$P_{max.c}(t) = V_{oc}(t)I_{sc}(t)FF(t) \qquad (1)$$

$$V_{oc.c}(t) = V_{oc.stc} + K_v(T(t) - 25) \qquad (2)$$

$$I_{sc.c}(t) = (I_{sc.stc} + K_i(T(t) - 25))\frac{G(t)}{1000} \qquad (3)$$

$$FF(t) = \frac{P_{max.c}(t)}{V_{oc.c}(t)I_{sc.c}(t)} \qquad (4)$$

here $V_{oc.c}(t)$ and $I_{sc.c}(t)$ are the calculated open circuit voltage and the short circuit current of the PV-cell, receptively. $V_{oc.stc}$ and $I_{sc.stc}$ are the given open circuit voltage and the short circuit current of the PV-cell under standard test conditions, receptively. $K_v$ is the given open circuit voltage temperature coefficient in the datasheet, $K_i$ is the given short circuit current temperature coefficient in the datasheet, $T$ is the actual PV-cell temperature coefficient and $FF(t)$ is the PV-cell fill factor. The detailed PV-cell model can be found in [12].

The instantaneous total available power from the whole PV-cells in a PV-array $P_{av.pv}(t)$ can be calculated by

$$P_{av.pv}(t) = N_{s.p}N_{p.p}N_{c.p}P_{max.c}(t) \qquad (5)$$

where $N_{s.p}$, $N_{p.p}$ are the quantity of series and parallel PV-panels connected in the PV-array, respectively, and $N_{c.p}$ is the quantity of PV-cells in one PV-panel. To prevent dispatching active power from the PV-array $P_{disp.pv}(t)$ more than the available power, the following constraint should be held

$$P_{disp.pv}(t) \leq P_{av.pv}(t). \qquad (6)$$

In this work, we assume that the PV-inverter is able to produce active as well reactive power. Therefore, the consumed AC power from the PV-inverter should be limited by PV-inverter rated apparent power as follows

$$\sqrt{P_{disp.pv}^2(t) + Q_{disp.pv}^2(t)} \leq S_{pv.inv} \qquad (7)$$

where $S_{pv.inv}$ is the PV-inverter apparent power rated capacity and $Q_{disp.pv}(t)$ is the reactive power dispatched from the PV-inverter .

### B. Battery bank

The battery bank is responsible to store the electrical energy to be used when it needed. The stored energy level in the battery bank is expressed by the state of charge ($SOC$) value which is increased by the charging power $P_{ch}(t)$ and decreased by the discharging power $P_{dis}(t)$ at each time step. The $SOC$ value is calculated as follows

$$SOC(t + \Delta t) = SOC(t) + s_1(t)\frac{\eta_{ch}P_{ch}(t)\Delta t}{E_{max}} \\ -s_2(t)\frac{P_{dis}(t)}{\eta_{dis}E_{max}}\Delta t \qquad (8)$$

here $E_{max}$ is the maximum energy of the installed battery bank. $s_1(t)$ and $s_2(t)$ are binary control variables used to prevent simultaneous charging and discharging through the

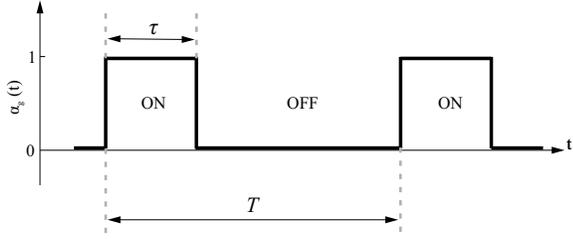

Fig. 2. The illustration of the grid status

battery bank operation. Therefore, the following constraints should be held

$$s_1(t) + s_2(t) \leq 1. \quad (9)$$

Moreover, the constraints on the $SOC$ should be defined to prevent overcharging or deep discharging, i.e.

$$SOC_{min} \leq SOC(t) \leq SOC_{max} \quad (10)$$

here $SOC_{min}$ is related to the battery bank depth of discharge $DOD$ by

$$SOC_{min} = (1-DOD)SOC_{max} \quad (11)$$

where

$$SOC_{max} = V_{b.n}Q_{b.n}. \quad (12)$$

Here $V_{b.n}$ and $Q_{b.n}$ are the nominal voltage and capacity in *ampere.hour* of the battery bank, respectively. In addition, the consumed AC power should be limited by the battery inverters apparent power capability [13]. Therefore, the following constraints are added

$$\sqrt{P_{ch}^2(t)s_1(t) + Q_{disp.b}^2(t)} \leq S_{b.inv} \quad (13)$$

$$\sqrt{P_{dis}^2(t)s_2(t) + Q_{disp.b}^2(t)} \leq S_{b.inv} \quad (14)$$

here $S_{b.inv}$ is the rated capacity of the battery inverter and $Q_{disp.b}(t)$ is the reactive power from the battery inverter.

*C. Diesel generator*

The fuel consumption of the DG $f_{con}(t)$ is related to the power dispatched from the DG $P_{disp.dg}(t)$ and the DG size [16], as follows

$$f_{con.dg_i}(t) = \begin{cases} AP_{disp.dg_i}(t) + BP_{n.dg_i}, & if \ P_{disp.dg_i}(t) > 0 \\ 0 & otherwise \end{cases} \quad (15)$$

where $A$ is a constant with the value of $0.246l/kWh$, $B$ is a constant with the value $0.08415l/kWh$. In Eq. (15), $P_{n.dg}$ is the nominal power of the DG.

Moreover, the reactive power capability of the DG is limited by manufacturer specifications to prevent overheating, therefore,

$$\sqrt{P_{disp.dg}^2(t) + Q_{disp.dg}^2(t)} \leq S_{dg.max} \quad (16)$$

$$PF_{r.dg} \leq PF_{dg}(t) \leq 1 \quad (17)$$

TABLE I
SIZE OF DEVELOPED MICROGRID COMPONENTS.

| Components | size |
|---|---|
| PV-modules | 700 kW˙p |
| PV-inverter | 700 kVA |
| Battery bank | 960 kWh |
| Battery inverter | 500 kVA |
| Diesel generator | 500 kVA |

where

$$PF_{dg}(t) = \frac{P_{disp.dg}(t)}{\sqrt{P_{disp.dg}^2(t) + Q_{disp.dg}^2(t)}} \quad (18)$$

where $S_{dg.max}$ and $PF_{r.dg}$ are the DG rated apparent power and power factor, respectively, $P_{disp.dg}(t)$ and $Q_{disp.dg}(t)$ are the dispatched active and reactive power from the DG, respectively, and $PF_{dg}(t)$ is DG power factor at time $t$.

*D. Grid blackouts*

In this paper, the studied microgrid is connected to a public grid that has a problem of scheduled blackouts [11], [12]. The available apparent power from the grid $S_{av.g}(t)$ is described as

$$S_{av.g}(t) = \alpha_g(t)S_{max.g} \quad (19)$$

where $S_{av.g}$ is the maximum available that can be imported from the public grid. $\alpha_g(t)$ is the grid availability parameter, when $\alpha_g(t) = 1$, the grid power is available and when $\alpha_g(t) = 0$, the grid power is unavailable. Here, the changing in the grid availability parameter is illustrated in Fig. 2. where $\tau$ is the Grid-ON period and $T$ is the total ON-OFF period.

It is to note that the following constraint should be held during importing active or/and reactive power from the grid

$$\sqrt{P_{disp.g}^2(t) + Q_{disp.g}^2(t)} \leq S_{g.max} \quad (20)$$

where $P_{disp.g}(t)$ and $Q_{disp.g}(t)$ are the active and reactive dispatched power from the public grid, respectively.

III. OPTIMIZATION PROBLEM FORMULATION

The objectives of the optimization are to minimize the cost of the consumed power by the load ($F_1$) and meanwhile to maximize the dispatched power from the PV-array ($F_2$). The optimization problem is defined as follows

$$\min_{\mathbf{u}_c(t), \mathbf{u}_b(t)} J = w_1 F_1 - w_2 F_2 \quad (21)$$

where

$$F_1 = C_{e.g} P_{disp.g}(t) + C_f f_{con}(t) \quad (22)$$

and

$$F_2 = P_{disp.pv}(t). \quad (23)$$

Here $\mathbf{u}_c(t)$ and $\mathbf{u}_b(t)$ are the continuous and binary control variables vectors, respectively. $w_i$ is a weighting factor, $C_{e.g}$ is the cost of the consumed energy from the grid in $\$/kWh$ and $C_f$ is the cost of consumed fuel by the DG in $\$/l$.

TABLE II
CONSUMED ENERGY FROM THE GRID AND DG IN TWO DAYS AT EACH SEASON.

| Seasons | Diesel only | | PV-battery-diesel | |
|---|---|---|---|---|
| | $E_{disp.g}(kWh)$ | $E_{disp.dg}(kWh)$ | $E_{disp.g}(kWh)$ | $E_{disp.dg}(kWh)$ |
| Winter | 5212.8 | 5212.8 | 3867.7 | 3687.5 |
| Spring | 4616.7 | 4616.7 | 2673.2 | 2801.7 |
| Summer | 7447.8 | 7447.8 | 3555.9 | 3705.5 |
| Fall | 4318.7 | 4318.7 | 2286 | 3163.2 |

TABLE III
COST OF THE CONSUMED ENERGY FOR EACH SCENARIO IN TWO DAYS.

| | Winter | Spring | Summer | Fall |
|---|---|---|---|---|
| Diesel only ($) | 4220.1 | 3910.8 | 5380.1 | 3756.1 |
| PV-battery-diesel ($) | 3013.7 | 2318.4 | 2847.4 | 2456.8 |
| Improvement | 28.59% | 40.72% | 47.08% | 34.59% |

IV. RESULTS AND DISCUSSIONS

To illustrate the potential of the developed MicroGrid model shown in Fig. 1 and the proposed optimal operation strategy in section III in solving the problem of grid blackouts. A case study is adapted from [12] with components sizes in Table I. The power factor of the connected industrial loads considered to be fixed and equal 0.85. The active power load profile is taken from [17] and shown in Fig. 3 (a). The blackout period is $\tau = 8$ hours and the total ON-OFF period $T = 16$ hours, see Fig. 2, meanwhile, the maximum capability of the grid connection is 500 kVA. The price of the dispatched active power form the grid is $C_{e.g}$ is 0.15 $/kWh. For the DG, only the fuel consumption cost is counted with 1.5 $/l. The computation results for a selected two days at each season of the same year are shown in Tables II - III and Fig. 3 (b) - (e).

As can be seen from Fig. 3, the proposed microgrid can provide active as well reactive power to the load while satisfying the microgrid technical constraints. From Fig. 3 (c), it can be noted that most output of the PV-array is used. From Fig. 3 (b) and (c) it can be seen thay the PV-array output power is used to cover the load even if the grid power is available. Also, the PV-array output power is used to decrease the dispatched power from the DG, see Fig. 3 (c) and (d). Moreover, the surplus power generated from the PV-array is stored in the battery bank which highly increases the total power consumeption from the PV-array, as shown in Fig. 3 (c) and (e). In addition, it is clean from Fig. 3 (d) and (e) that when the load is lower than the minimum load limit of the DG [18], the battery bank works as a dump load to increase the consumed power from the DG.

As a result of taking into account the consumed energy from the grid and the DG costs in the objective function of the optimization problem, a significant reduction of the consumed energy from the grid $E_{disp.g}$ and the DG $E_{disp.dg}$ is gained, as given in Table II. The costs of the consumed energy from the grid and the DG in details are given in table III. The obtained results indicate that using the proposed microgrid instead of a stand-alone DG reduces the power consumption cost up to 47.08%.

V. CONCLUSION

Microgrids are promising solutions to guarantee a reliable electrical power source for areas suffering from grid blackouts. In this work, a new model for a PV-battery-diesel microgrid is proposed to provide a reliable electrical power source for industrial loads considering both active and reactive power flow. Moreover, an optimal operation strategy is introduced to decrease the consumed energy cost and increase the dispatched energy from the PV-array. The result show that using the proposed microgrid can provide an uninterrupted power supply to the load and considerably reduce the total dispatched energy cost. In addition, the proposed optimal operation strategy can reduce the power dispatched from the grid and the DG while decreasing the power curtailment from the PV-array, that maximizes the economic and environmental benefits of installing the PV-array. On the basis of the promising findings presented in this paper, the effect of the consumed reactive power cost on the operation of the microgrid will be involved in future work.

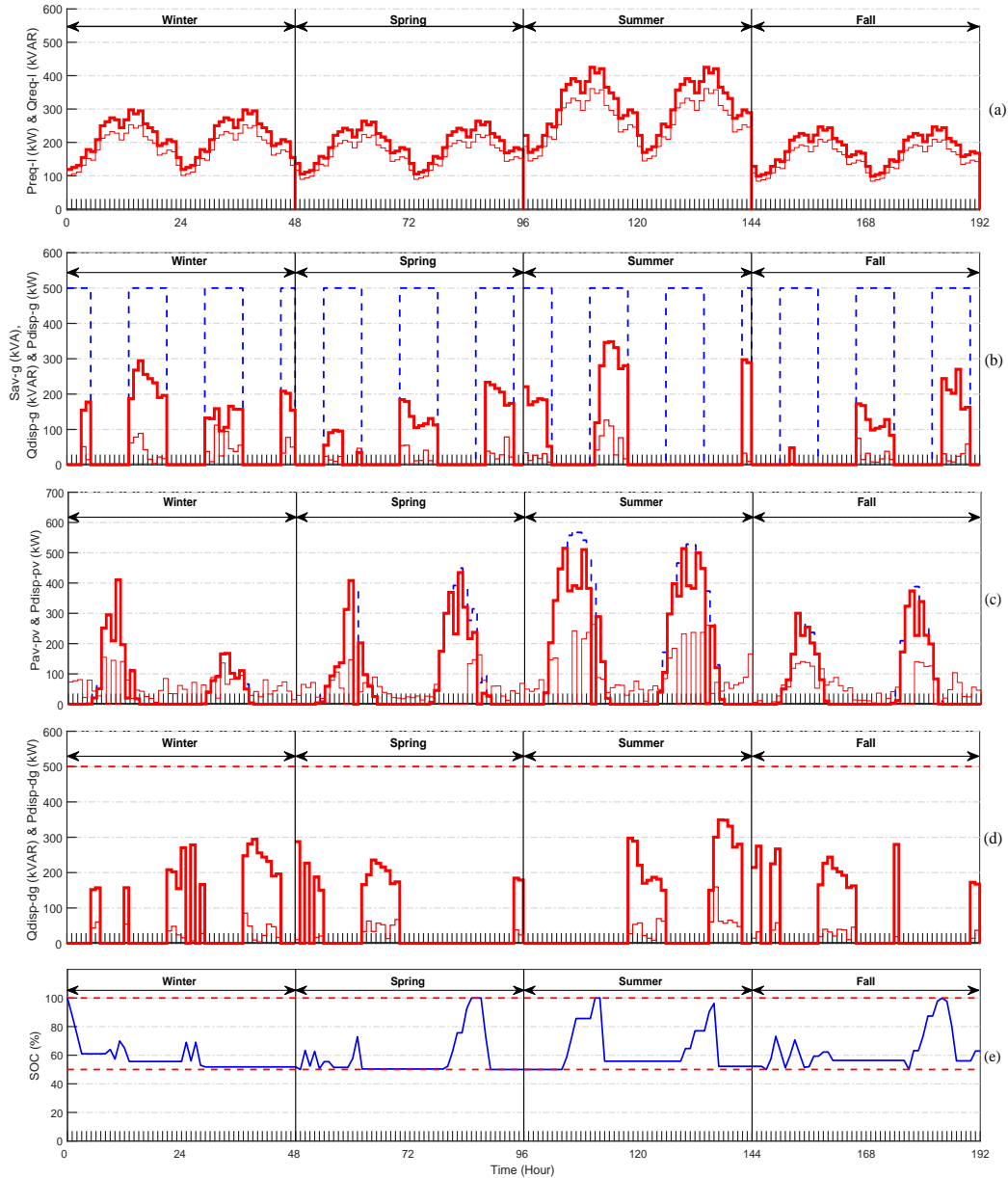

Fig. 3. Microgrid components behavior in two days at each season . (a) Load active power (bold-red) and load reactive power (thin-red). (b) Available power (dashed-blue), active power dispatch (bold-red), reactive power dispatch (thin-red) from grid. (c) Available power (dashed-blue) active power dispatch (bold-red), reactive power dispatch (thin-red) from PV-array. (d) Active power dispatch (bold-red), reactive power dispatch (thin-red) from DG. (e) Battery bank state of charge (solid-blue) and its upper and lower limits (dashed-red).